\documentclass[10pt]{article}

\usepackage{graphicx}			
\usepackage{amssymb}
\usepackage{amsmath}
\usepackage{psfrag}
\usepackage{verbatim}
\usepackage{xcolor}

\parskip = \baselineskip

\def\vp{\vspace{\baselineskip}}

\title{The Trilobite and Crab, a full explanation}
\date{}
\author{Chaim Goodman-Strauss\\ Univ. Arkansas\\ \tt strauss@uark.edu}

\begin{document}

\maketitle

\centerline{\includegraphics[]{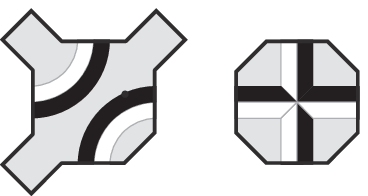}}

The ``trilobite and crab", shown above, admit tilings of the plane but admit {\em only} non-periodic tilings, and so are an {\em aperiodic} set of tiles. They are 
 among the very simplest aperiodic set of aperiodic tiles known in ${\mathbb E}^2$  --- there are only a few other pairs known and only one other,  Amman's A2~\cite{grsh},   has as few translation classes, i.e. appears in so few orientations, eight.  (It is still a well-known open question whether or not there is a single aperiodic tile.)
 A complete  bibliography appears in~\cite{gs_LASTs}.

The pair is derived from the ``trilobite and cross" tiles, described in~\cite{gs_small} (which  generalize to an  aperiodic pair of  tiles in all ${\mathbb E}^{n\geq 3}$). The proof that the trilobite and cross tiles are aperiodic is a fairly simple combinatorial check that the tiles {\em can} form larger patches with the same combinatorial structure, which can then be assembled into still larger patches, {\em ad infinitum}, and thus can tile the plane. 
Conversely, in any tiling by these tiles, the trilobites {\em must} lie in such a hierarchy of patches, showing no such tiling can have a translational period. 

However, like the ``Pegasus" pair of tiles~\cite{gs_pegasus},  the trilobite and cross have unusual ``tip-to-tip" matching rules; we can easily recompose these into three tiles with matching rules that are completely encoded geometrically (by ``bumps and nicks")  or by colored edges, as shown at right  below. (The tiles can be adjusted to have areas 1, $\epsilon$ and $\epsilon^2$, as in~\cite{penroseEpsilon}.)

\centerline{\includegraphics[]{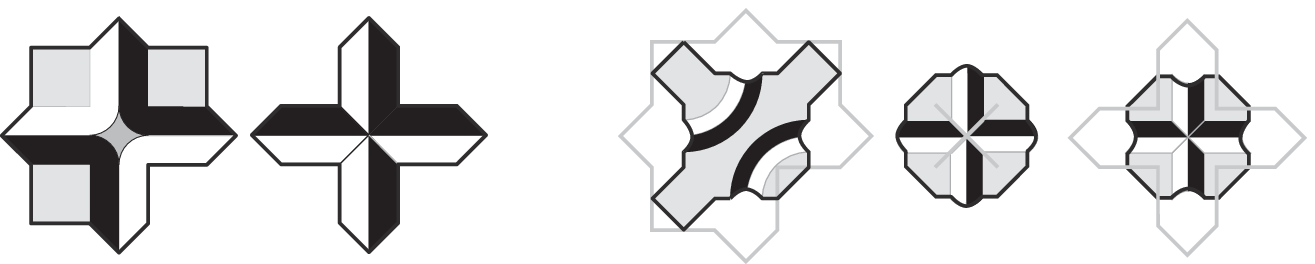}}

This raises the natural question: can we conflate two of these tiles, giving us the crab tile, and  have an aperiodic pair? We certainly will allow a richer variety of local configurations, giving more complex combinatorial structure.

It is {\em quite remarkable} that the proof complexity of the following theorem seems to be quite high:

{\bf Theorem} {\em The trilobite and crab are an aperiodic pair of tiles.}

Because of the undecidability of the domino problem, it is certainly the case that as we enumerate all possible sets of tiles, among those that are aperiodic, the length of the shortest  proof that they are so cannot be bounded by any computable function (see~\cite{gs_cantdecide}). But is amazing, at least to me, that this kicks in so readily.  

As I wrote in~\cite{gs_small}, a full proof of this theorem is ``not worth the readers time", but it is worth having as a striking example of this phenomenon. Moreover, there seem to be many interesting possibilities for exploiting this complexity, such as programming within defects of the tiling. 

So, at last, here is a full proof, given in a graphical shorthand, drawn on a square grid. 

{\em Proof of the Theorem:} 

We will encode the trilobite and crab thusly: 

\centerline{\includegraphics[]{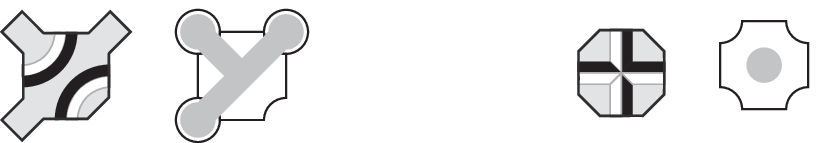}}

Here is an encoding of a typical configuration:

\centerline{\includegraphics[]{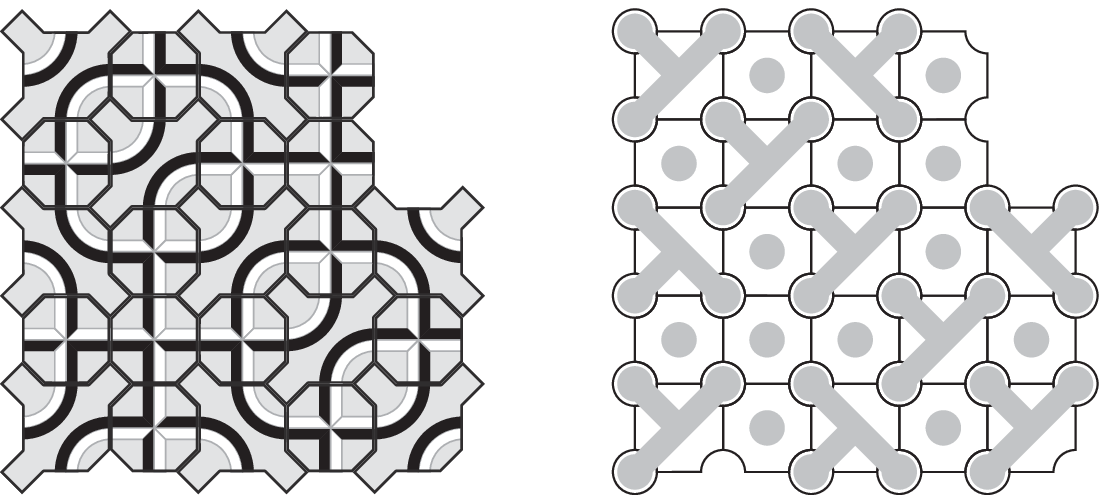}}

Our proof is essentially that the matching rules enforce hierarchical configurations such as this one:

\centerline{\includegraphics[]{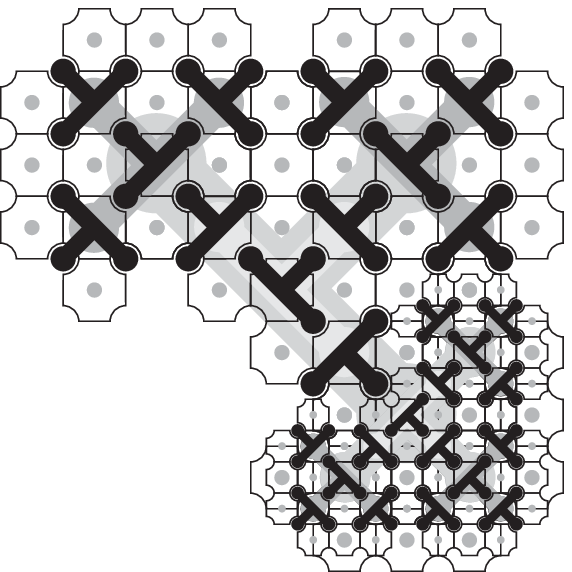}}

 We adopt several further conventions: We take  \includegraphics[]{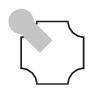} to mean \includegraphics{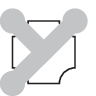}, \includegraphics{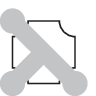}, or \includegraphics{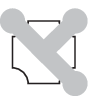}.  We take \includegraphics[]{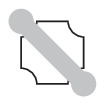} to mean \includegraphics{pix/trilobiteAndCrab_code4c.eps}, or \includegraphics{pix/trilobiteAndCrab_code4d.eps}, and we take \includegraphics{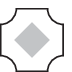} to mean a trilobite in any orientation.
 
In each case of the proof, solid {\color{black} {\bf black}} objects are given, and {\color{gray} {\bf gray}} ones are implied, with numbers giving the order of the implication; {\color{green}{\bf green}} indicates reduction to an earlier case, {\color{blue} {\bf blue}} means a further subcase, and {\color{red} {\bf red}} indicates a contradiction. 
 
 \paragraph{Axioms}  With \includegraphics{pix/trilobiteAndCrab_code4b.eps} and \includegraphics{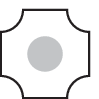}, \begin{enumerate} \item Each square and corner of a square grid must be covered, and 
 \item By parity, any pair of trilobites separated by crabs must be oriented as: {\includegraphics{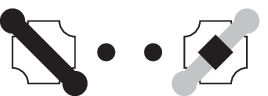}} \end{enumerate}
 
 We now have 
 
 \paragraph{Several Elementary Lemmas} (The helpful paper included at the end of these notes will be useful for checking these and all our arguments.)
 
{\includegraphics[]{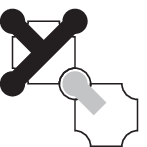}}, meaning that adjacent to   {\includegraphics[]{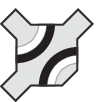}}, we must have  {\includegraphics[]{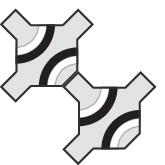}},  {\includegraphics[]{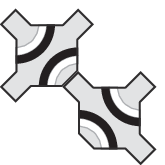}}, or  {\includegraphics[]{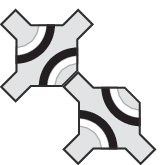}}.

Similarly, each of:

\centerline{\includegraphics[]{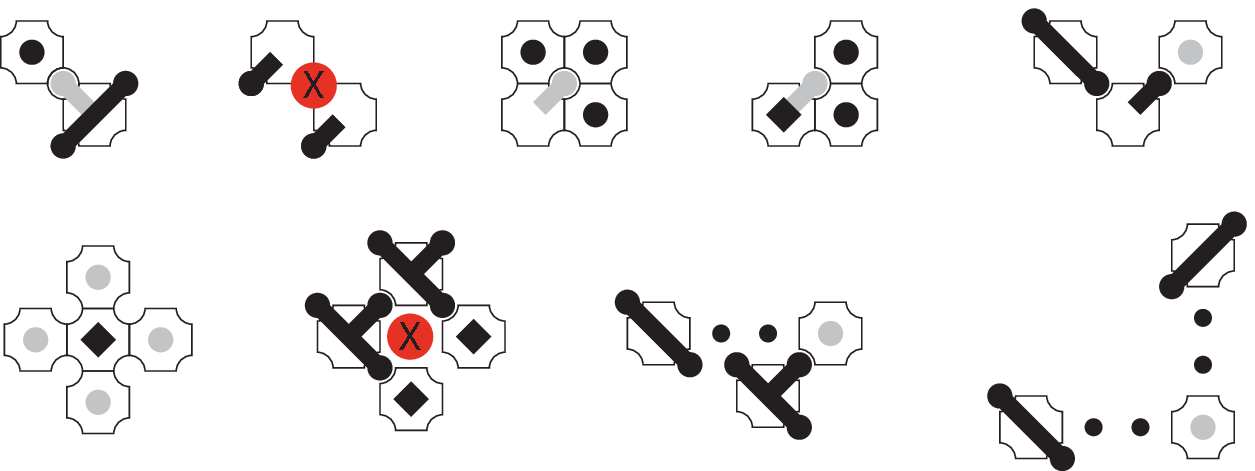}}

\paragraph{Enumerating initial cases} We consider the possible tiles around a trilobite, naming the configurations similarly to~\cite{gs_small}, writing {\tt T} for trilobite, {\tt O} for crab, or {\tt *} for either.

\psfrag{a}{\tt ***}\psfrag{b}{\tt T**}\psfrag{c}{\tt O**}\psfrag{d}{\tt *T*}\psfrag{e}{\tt *O*}

\centerline{\includegraphics[]{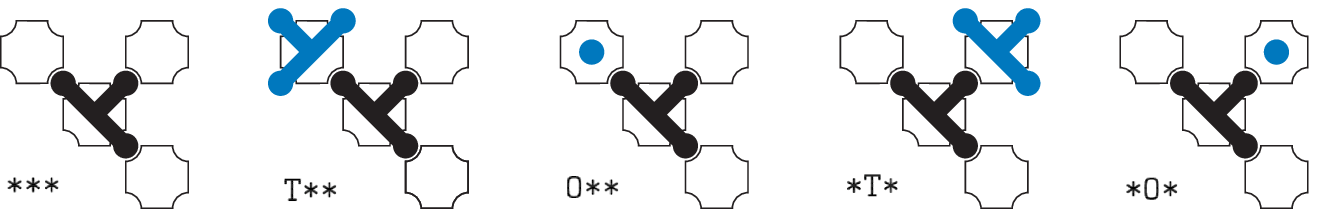}}

\vp {\tt TTT}, {\tt OTO} and {\tt OOO} arise within the combinatorial structure we seek:

\centerline{\includegraphics[]{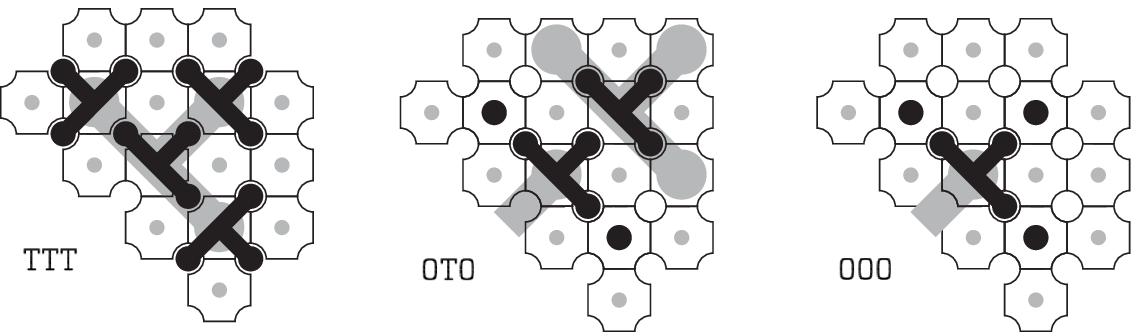}}

However (as in~\cite{gs_small}) we must take special care with {\tt OTT} and {\tt TTO} and ensure that   {\tt OOT}, {\tt TOT} and {\tt TOO} are forbidden entirely.

\paragraph{{\tt TOT} is forbidden:}
{\includegraphics{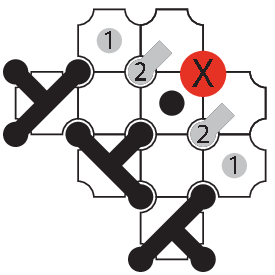}}

\paragraph{{\tt *OT} is forbidden:} {\includegraphics[]{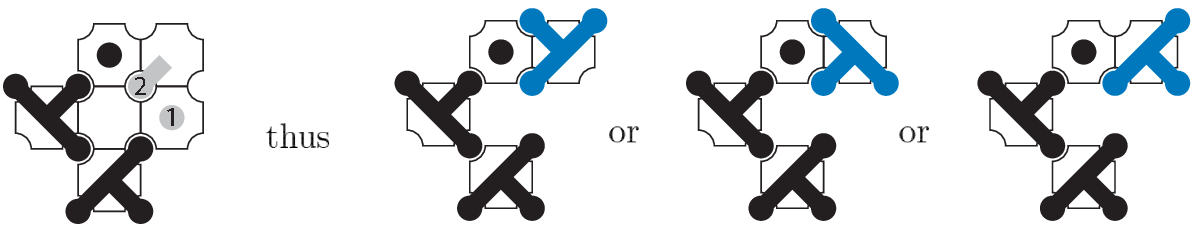}}

Taking each of these cases in turn, and noting that that the later cases reduce to the first (but indicating all of the implied tiles), we have:

\centerline{\includegraphics[]{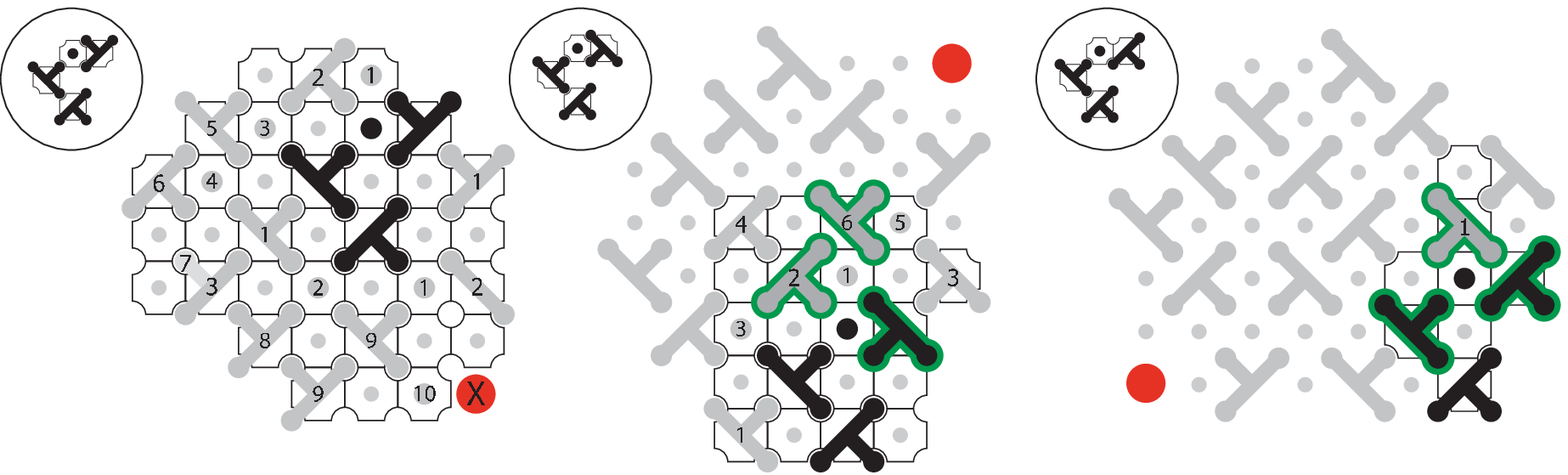}}

\paragraph{Chains of {\tt *TO}'s are forbidden}

As in~\cite{gs_small}, the trilobite and crab do admit chains of alternating {\tt TTO}'s and {\tt OTT}'s, but we must show that any {\tt *TO} can only appear in this way.

\centerline{\includegraphics[]{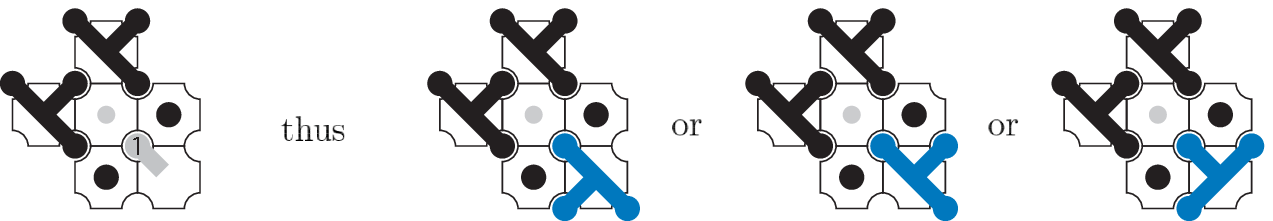}}

\centerline{\includegraphics[]{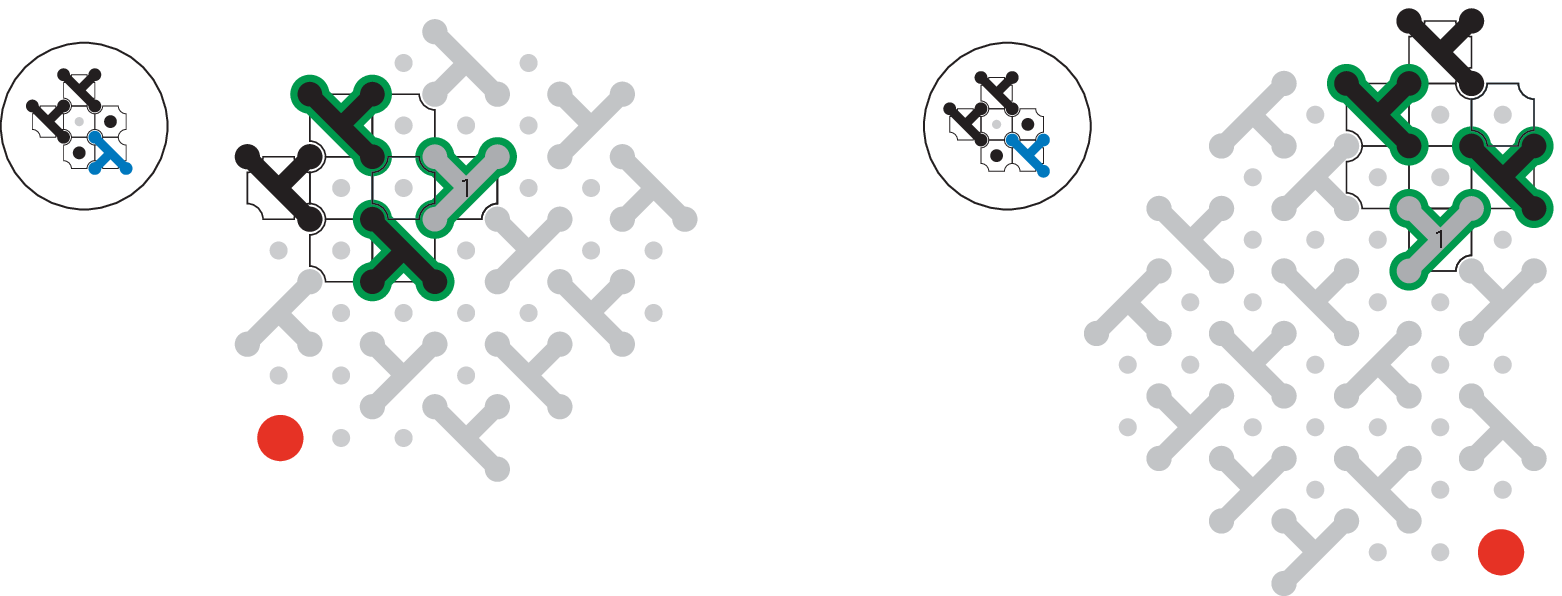}}

This case is more complex still: {\includegraphics[]{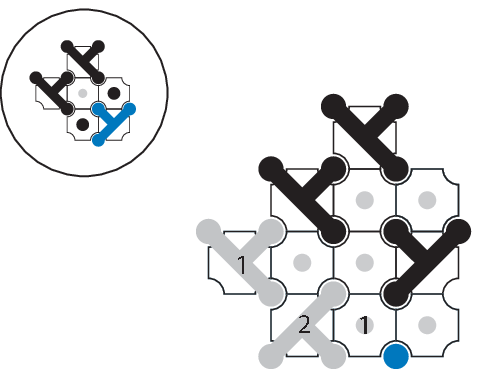}}

{\includegraphics[]{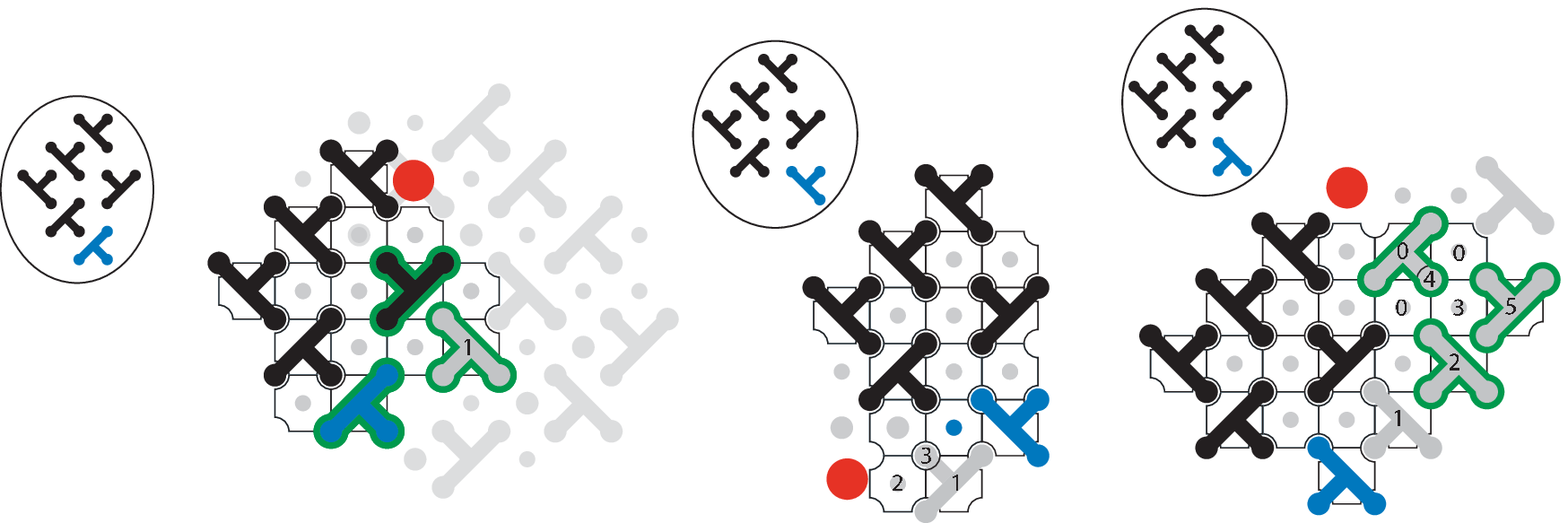}}

Finally we can conclude:

\centerline{\includegraphics[]{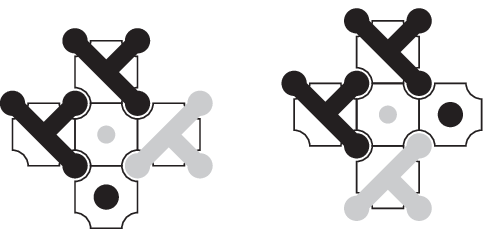}}

Any {\tt *TO} in a tiling must therefore occur in an infinite chain of alternating {\tt TTO}'s and {\tt OTT}'s and as in~\cite{gs_small}, there must always be a corresponding tiling, formed by sliding half the plane one tile along this diagonal. In this new tiling, all trilobites will be {\tt TTT}, {\tt OTO} or {\tt OOO}.

\centerline{\includegraphics[]{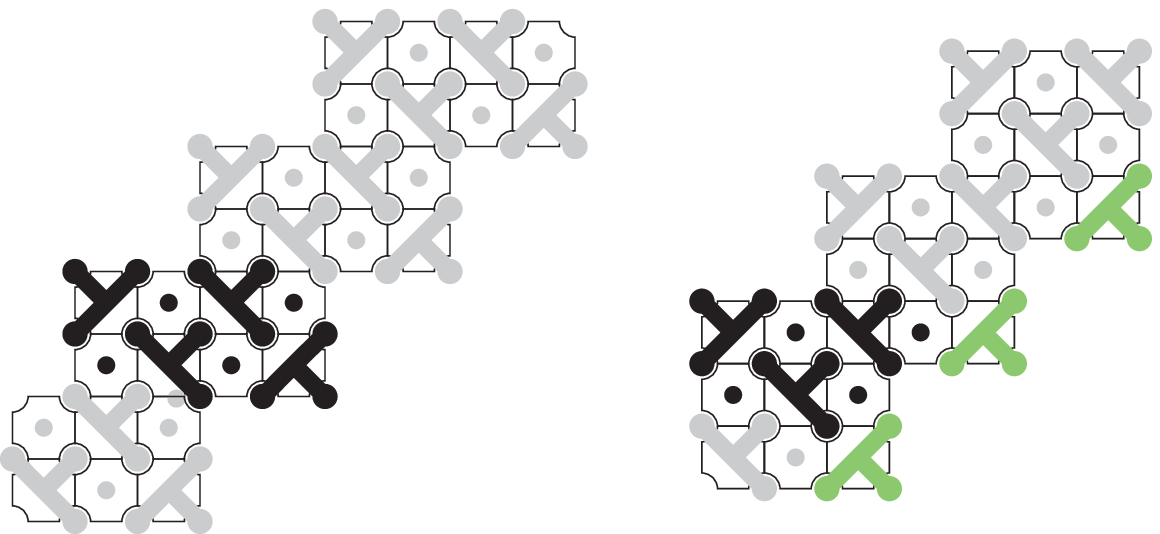}}

\paragraph{Establishing the induction}
Every trilobite of type {\tt TTT} in effect is a larger trilobite, but we {\em still} must check that these large tiles satisfy our axioms. The difficulty is ensuring that only the top left alignment is allowed:

\centerline{\includegraphics[]{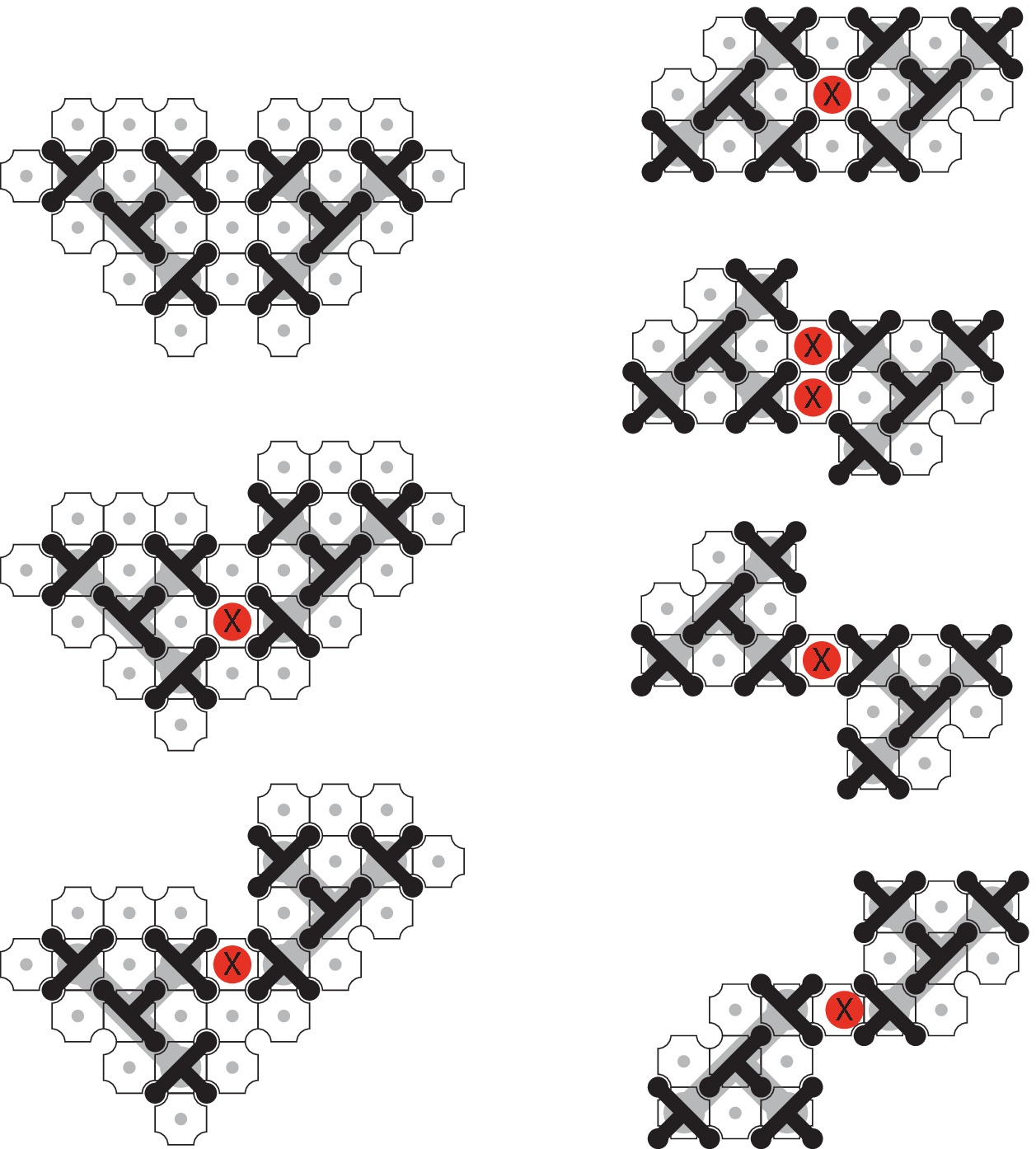}}

One last case breaks  into still further subcases:

\centerline{\includegraphics[]{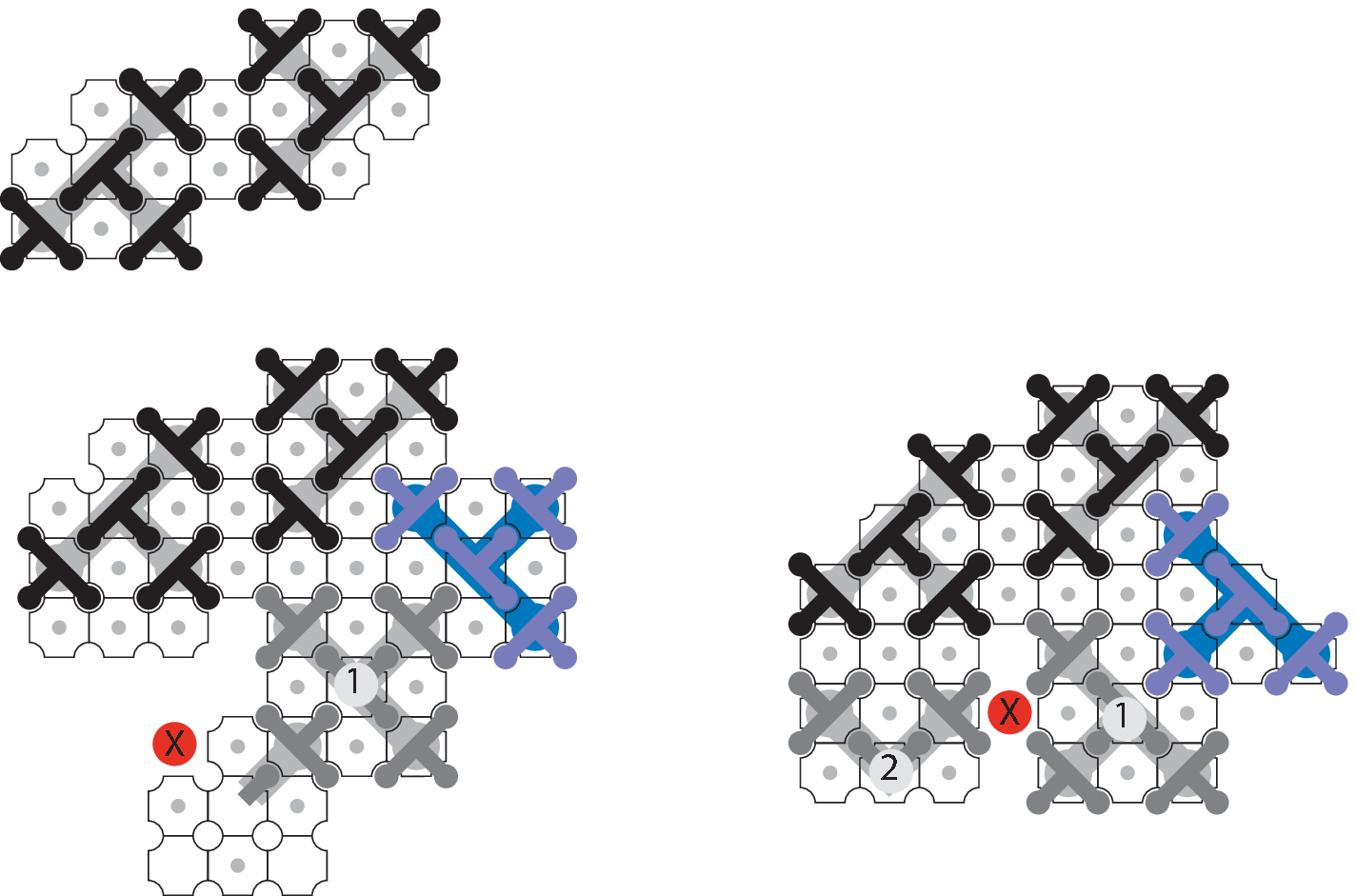}}

Finally, after all of this, we have that, in any tiling by the trilobite and crab, every trilobite is in a larger trilobite, or after a shift along a chain of {\tt *TO}'s this is so, and that these larger trilobites satisfy the axioms as before. 

\centerline{\includegraphics[]{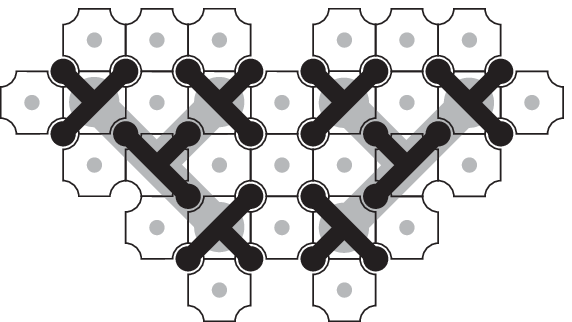}}

We may therefore  induct, and the proof, as in~\cite{gs_small}, is {\em finally} complete!

\centerline{\includegraphics[]{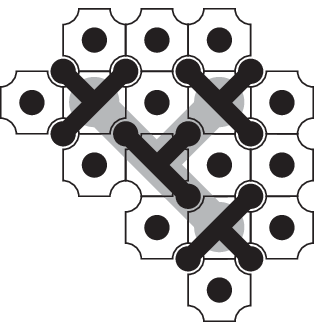}}

\centerline{\includegraphics[]{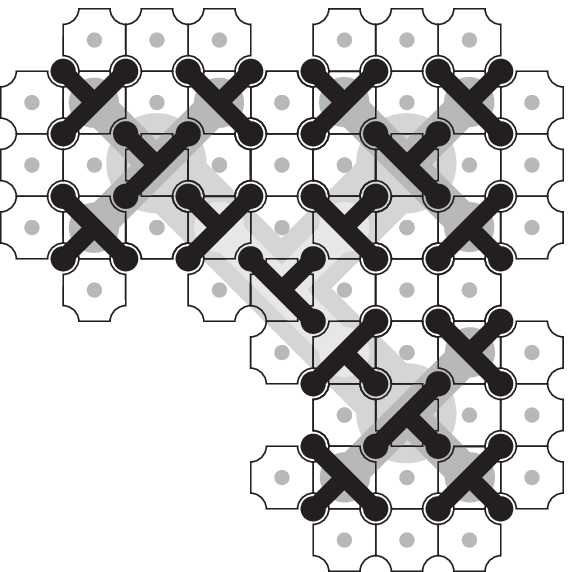}}

\centerline{\includegraphics[width=1.1\textwidth]{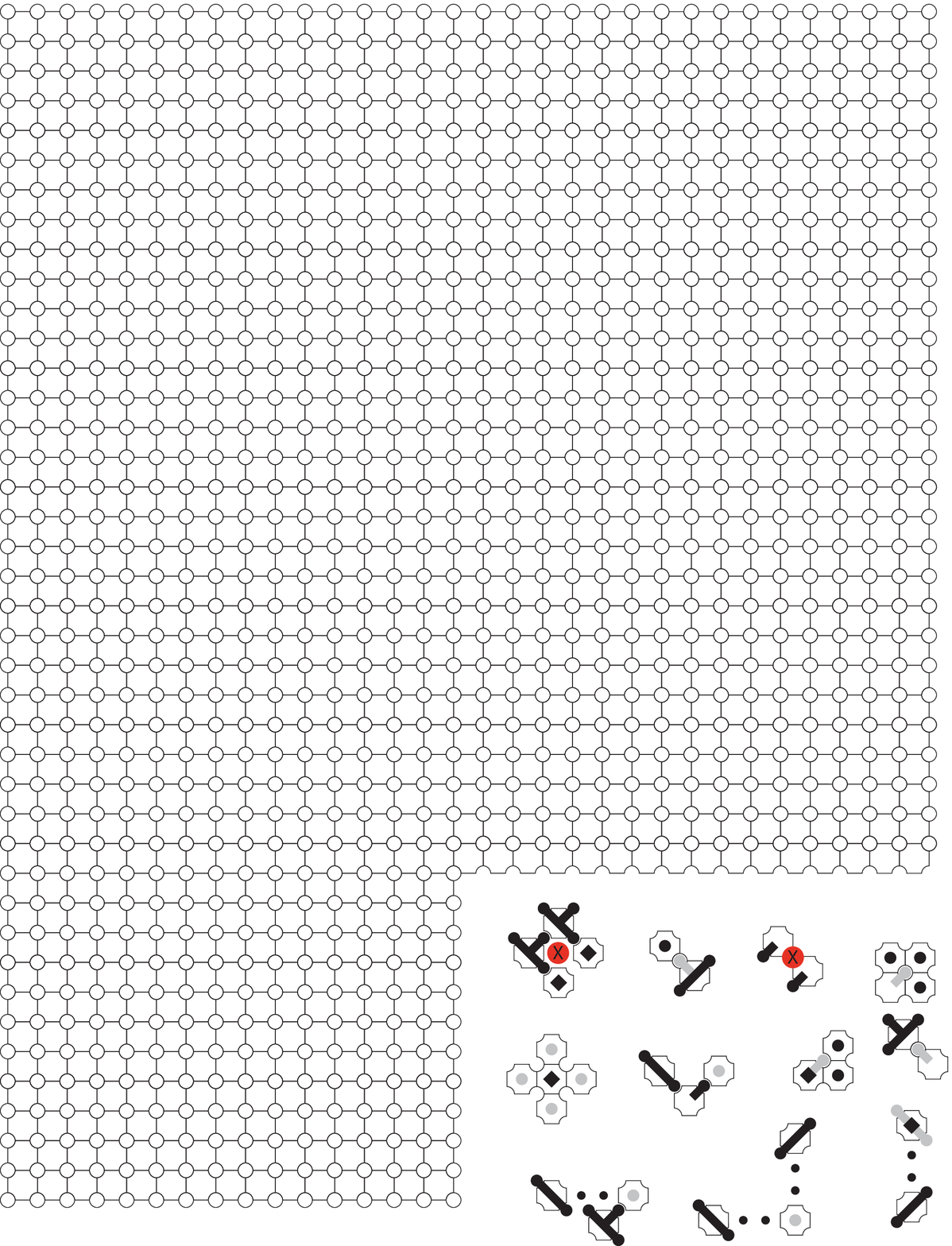}}

\centerline{\bf Useful paper for checking cases}

\small
\begin{thebibliography}{6}

\bibitem{ammanGrsh} {R.  Amman, B.  Grunbaum and G.C.  Shepherd}, 
{\em Aperiodic tiles}, Discrete and Computational Geometry {\bf 8} 
(1992) 1-25.

\bibitem{berger} R. Berger, {\em The undecidability of the domino problem},
Memoirs Am. Math. Soc.  {\bf 66} (1966).


\bibitem{gs_mrst} {C. Goodman-Strauss}, {\em Matching rules and substitution
tilings},  Annals of Math. {\bf 147} (1998), 181-223.


\bibitem{gs_small} {C. Goodman-Strauss}, {\em A small  aperiodic 
set of planar  tiles},  Europ. J. Combinatorics {\bf 20} (1999), 
375-384.

\bibitem{gs_en} {C. Goodman-Strauss}, {\em An aperiodic pair of 
tiles in  $E^n$
for all
$n\geq 3$},  Europ. J. Combinatorics {\bf 20} (1999), 385-395.

\bibitem{gs_cantdecide} C. Goodman-Strauss, {\em Can't Decide? Undecide!}, Notices A.M.S. {\bf 57} (2010), 343-356. 

  \bibitem{gs_LASTs}{C. Goodman-Strauss}, {\em Lots of Aperiodic Sets of Tiles}, arXiv

  \bibitem{gs_pegasus}{C. Goodman-Strauss}, {\em The Pegasus tiles}, arXiv
  
 \bibitem{gs_fullExplanation}{C. Goodman-Strauss}, {\em Matching rules for the sphinx substitution tiling}, arXiv


\bibitem{grsh}{B. Gr\"unbaum} and {G.C. Shepherd}, {\em Tilings and patterns},
W.H. Freeman and Co.   (1987).


 \bibitem{penroseEpsilon} R. Penrose {\em Remarks on Tiling: details of a
$(1+\epsilon+\epsilon^2)$-aperiodic set}, The mathematics long range
aperiodic order, NATO Adv. Sci. Inst. Ser. C. Math. Phys. Sci. 489 (1997),
467-497.



\end{thebibliography}
\end{document}